%dyck.tex: Automatic Counting of Restricted Dyck Paths via (Numeric and Symbolic) Dynamic Programming 
%%a Plain TeX file by Shalosh B. Ekhad and Doron Zeilberger (x pages)

%begin macros

\baselineskip=14pt
\parskip=10pt

\magnification=\magstephalf

\def\P{{\cal P}}

\def\1{{\overline{1}}}
\def\2{{\overline{2}}}
\parindent=0pt
\overfullrule=0in

\def\frac#1#2{{#1 \over #2}}
%\headline={\rm  \ifodd\pageno  \RightHead  \else  \LeftHead  \fi}
%\def\RightHead{\centerline{
%Title
%}}
%\def\LeftHead{ \centerline{Doron Zeilberger}}
%end macros
\centerline
{\bf  
Automatic Counting of Restricted Dyck Paths 
}
\centerline
{\bf via (Numeric and Symbolic) Dynamic Programming }

\bigskip
\centerline
{\it Shalosh B. EKHAD and Doron ZEILBERGER}
\bigskip
\quad\quad\quad\quad {\it In memory of our hero Richard Guy (1916-2020)}
\bigskip
\quad\quad {\it ``On some fields it is difficult to tell whether they are sound or phony. Perhaps they
are both. Perhaps the decision depends on the circumstances, and it changes with time $\dots$. Some subjects
 start out with impeccable credentials, catastrophe theory for instance, and then turn out to resemble
 a three-dollar bill. Others, like dynamic programming, have to overcome a questionable background before they
are reluctantly recognized to be substantial and useful.''} \hfill\break
-- Gian-Carlo Rota (``Discrete Thoughts'', Birkh\"auser, 1992, p. 1)

{\bf Preface} 

In his classic essays ``The Strong Law of Small Numbers" [G1][G2], Richard Guy gave numerous
{\it cautionary tales} where one can't `jump to conclusions' from the first few terms of a sequence.
But if you are cautious enough you can find many families of enumeration problems where
it is very safe to deduce the general pattern from the first few cases, obviating the need
for either the human or the computer to think too hard, and by using the {\bf `Keep It Simple Stupid'} 
principle ({\bf KISS} for short), one can easily derive many deep enumeration theorems by doing
exactly what Richard Guy told us {\bf not} to do, compute the first few terms of the sequence
and deduce the formula for the general term. We admit that often  `few' should be replaced by
`quite a few', but it is still much less painful than trying to figure out the intricate 
combinatorial structure by `conceptual' means.

The way we generate sufficiently many terms of the studied sequence that enables guessing the pattern
is via {\it numeric} {\bf dynamic programming}. While it is very fast to generate many terms
of the sequence, the guessing part, for larger problems, would eventually  become impractical.
For these larger problems, both the human programmer and the machine need to think a bit,
and use {\bf symbolic} dynamic programming. Note that the human only has to think {\bf once},
writing a general program that teaches the computer how to `think' in each specific case of
a large class of enumeration problems.

{\bf The Maple packages} 

This article is accompanied by two Maple packages

$\bullet$ {\tt Dyck.txt}, that uses {\bf Numeric} Dynamic programming 
and {\it number-crunching}
to generate many terms
of the desired sequence, and then uses the KISS method to guess algebraic equations and
linear recurrences that can be rigorously justified {\it a posteriori}.

$\bullet$ {\tt DyckClever.txt}, that uses {\bf Symbolic} Dynamic programming 
and {\it symbol-crunching} to directly get the
desired equations for the generating functions of interest, and by using further symbol-crunching,
derives linear recurrences for the sequences of interest.

Both packages are available, along with input and output files with ample data, from the front of the present article

{\tt http://www.math.rutgers.edu/\~{}zeilberg/mamarim/mamarimhtml/dyck.html} \quad .

{\bf Part I: Using Numeric Dynamic Programming to Enumerate Restricted Dyck Paths}

{\bf Enumerating Dyck Paths}

Recall that a {\bf Dyck path} of {\bf semi-length} $n$ is a walk in the plane,  from the origin
$(0,0)$ to $(2n,0)$ with {\bf atomic steps} $U:=(1,1)$ and $D:=(1,-1)$ that {\bf never} goes below the
$x-$axis, i.e. that always stays in $y \geq 0$.

For example, the five Dyck paths of semi-length $3$ are
$$
UUUDDD \quad, \quad
UUDUDD \quad, \quad
UUDDUD \quad, \quad
UDUUDD \quad, \quad
UDUDUD \quad.
$$

The number of Dyck paths of semi-length $n$ is famously the {\bf Catalan number}, $\frac{(2n)!}{n!(n+1)!}$,
the most popular, and {\bf important} sequence in enumerative combinatorics (with no offense to Fibonacci),
the subject of a whole book by Guru Richard Stanley [St].

When we searched (on June 2, 2020) our favorite website, the OEIS [Sl] for the phrase ``Dyck paths" we got
back 1100 hits. Of course we did not have the patience to read  all of them, but a random browsing
revealed that they enumerate Dyck paths with various {\bf restrictions}, and they refered to
the interesting papers [ELY], [PW], and [M].
These were all (admittedly clever) human-generated efforts.
We will soon see how to quickly {\it automatically} enumerate `infinitely' many such classes, but first let's recall
one of the many proofs of the fact that the number of Dyck paths of semi-length $n$ is indeed the venerable OEIS sequence A108,
$\frac{(2n)!}{n!(n+1)!}$.

Let $a(n)$ be the desired number, i.e. the number of Dyck paths of semi-length $n$, and consider the
{\bf ordinary generating function}
$$
f(x) \, := \, \sum_{n=0}^{\infty} \, a(n) \, x^n \quad ,
$$
which is the {\bf weight-enumerator} of the set of all Dyck paths, with $weight(P):=x^{SemiLength(P)}$.

It is readily seen that any Dyck path $P$ is either empty or can be written {\bf uniquely}
(i.e. {\bf unambiguously}) as $P=U\,P_1\,D\,P_2$, where
$P_1$ and $P_2$ are shorter Dyck paths, and vica versa,  for any Dyck paths $P_1$, $P_2$,  $U\,P_1\,D\,P_2$ is
a Dyck path on its own right. Let $\P$ be the {\bf totality} of all Dyck paths, then we have the
{\bf grammar}
$$
\P \, = \, \{ EmptyPath \} \, \cup \, U\, \P \, D \, \P \quad.
$$
Applying the {\bf weight} functional we get
$$
f(x)= 1 + xf(x)^2 \quad .
$$

To deduce that $a(n)=(2n)!/(n!(n+1)!)$,  you can, {\it inter alia}  

$\bullet$ (i) Solve the quadratic and use Newton's binomial theorem  .

$\bullet$ (ii) Differentiate both sides getting a differential equation for $f(x)$ that translates to a first-order recurrence for $a(n)$.

$\bullet$ (iii)  Use Lagrange Inversion (see [Z1] for a brief and lucid  exposition).

{\bf How it All Started:  Vladimir Retakh's Question}

Volodia Retakh asked whether there is a proof of the fact that the number of Dyck paths
of semi-length $n$ such that the height of all peaks is either $1$ or even is given
by the also famous Motzkin numbers (OEIS sequence A1006), whose generating function satisfies
$$
f(x)=1 \,+ \, x f(x) \,+ \, x^2 f(x)^2 \quad .
$$
We first tried to find a `conceptual' proof  generalizing the above proof, and indeed we found one,
by adapting the above classical proof enumerating all Dyck paths.

Let $\P_1$ be the set of Dyck paths whose peak-heights  are never in $\{ 3,5,7, \dots \}$,
and let $f_1=f_1(x)$ be its weight enumerator.

Let $P_1$ be any member of $\P_1$ then,  it is either empty, or we can write
$$
P_1 \, = \, U \, P_2 \, D \, P'_1 \quad,
$$
where $P'_1 \in \P_1$ but $P_2$  has the property that none of its peak-heights is in $\{2,4,6, \dots \}$.
Let $\P_2$ be the set of such Dyck paths.

Hence, we can write the `grammar'
$$
\P_1 \, = \, \{ EmptyPath \}\, \cup \, U\, \P_2 \, D \, \P_1 \quad.
$$

Let $f_2=f_2(x)$ be the weight-enumerator of $\P_2$.

Taking weights above, we have the equation
$$
f_1 =  1 + x f_2 f_1 \quad .
$$

Alas, now we have to put-up with $\P_2$ and $f_2$. Let $P_2$ be any member of $\P_2$. Then
either it is empty, or  it can be written as
$$
 U \, P_3 \, D \, P'_2 \quad,
$$
where $P'_2 \in \P_2$ but $P_3$ is a Dyck path whose peak-heights are  never in $\{1,3,5,7,\dots\}$.

Let $\P_3$ be the set of such Dyck paths. We have the grammar
$$
\P_2 \, = \, \{ EmptyPath \}\, \cup \, U\, \P_3 \, D \, \P_2 \quad.
$$

Let $f_3=f_3(x)$ be the weight-enumerator of $\P_3$.

Taking weights we have another equation
$$
f_2 = 1 + xf_3 f_2 \quad .
$$

It looks like we are doomed to have {\bf infinite regress}, but let's try one more time.

Let $P_3$ be any member of $\P_3$. It is either empty, or we can write
$$
P_3 \, = \, U \,P_4\,D P'_3 \quad
$$
where $P'_3 \in \P_3$ and $P_4$ is a {\bf non-empty} path avoiding peak-heights in $\{2,4,6,8, \dots \}$. But this
looks familiar, so the set of $P_4$-paths is really $\P_2 \backslash \{EmptyPath\}$, and we have the
grammar
$$
\P_3 \, = \, \{ EmptyPath \} \, \cup \, U\, (\P_2 \backslash \{EmptyPath\}) \, D \, \P_3 \quad.
$$
Taking weight, we get
$$
f_3 = 1 + x(f_2-1) f_3 \quad .
$$

We have a system of three algebraic equations
$$
\{ 
f_1 = 1+xf_2 f_1 \quad, \quad
f_2 = 1+xf_3 f_2 \quad, \quad
f_3 = 1+x (f_2-1)f_3 \quad \} \quad,
$$
in the unknowns
$$
\{ f_1 , f_2 , f_3 \} \quad .
$$
Eliminating $f_2,f_3$ yields the following algebraic equation for our object of desire $f_1$,
$$
f_1(x)=1 \,+ \, x f_1(x) \,+ \, x^2 f_1(x)^2 \quad ,
$$
proving Volodia Retakh's claim.

In Part II we will see how to teach the computer to do these reasonings, but if neither human nor machine
feel like thinking too hard, for many problems one can use the {\bf KISS} method.

{\bf The KISS way}

Now that we know that such an argument {\bf exists}, and that the desired generating function $f(x)$, satisfies an algebraic
equation of the form $P(x,f(x))=0$ for {\it some} bivariate polynomial $P(x,y)$, why not {\bf keep it simple}, and
rather than wrecking our brains (both those of humans and those of  machines), let's collect sufficiently many terms of the
desired sequence, and then use Maple's command {\tt gfun[listtoalgeq]}  (or our own home-made
version) to {\bf guess} the desired  polynomial equation  $P(x,f(x)) \, = \, 0$. 

\vfill\eject

{\bf Numerical Dynamic Programming to the rescue}

Suppose that we don't know anything, and want to compute the number of Dyck paths of semi-length $n$, i.e. the
number of {\bf all} walks using the fundamental steps $U=(1,1)$ and $D=(1,-1)$, that start at $(0,0)$, end at $(2n,0)$
and never visit $y<0$.
A natural approach is to consider the more general quantity $d(m,k)$, the number of walks from $(0,0)$ to $(m,k)$ staying
weakly above the $x$-axis and {\bf ending at a down step}. 
If the length of that downward run is $r$, then the previous peak was at $(m-r,k+r)$, and
we need to introduce the auxiliary function $u(m,k)$ the number of such paths that end at $(m,k)$ and end at an up step.

We have
$$
d(m,k) \, = \,\sum_{r=1}^{m} u(m-r,k+r) \quad .
$$
Analogously,
$$
u(m,k) \, = \, \sum_{r=1}^{m} d(m-r,k-r) \quad .
$$

Of course we have the obvious {\bf initial condition} $d(0,0)=1$, and the {\bf boundary conditions} 
$d(m,k)=0$ and $u(m,k)=0$ if $k>m$.

Here is the short Maple code that does it

{\tt u:=proc(m,k) local r: option remember: if m=0 then  0: else add(d(m-r,k-r),r=1..k): fi: end:}

{\tt d:=proc(m,k) local r: option remember: \hfill\break
if m=0 then  if k=0 then    RETURN(1):  else   RETURN(0): fi: fi: add(u(m-r,k+r),r=1..m): end:
}

To get the desired sequence enumerating all Dyck paths of semi-length $n$ for $n$ from $1$ to $N$, in other words
$\{d(2n,0)\}_{n=0}^{N}$ for any desired $N$,  we type

{\tt seq(d(2*n,0),n=0..N);}   \quad .

Now, recall that we had to work much harder, {\it logically} and {\it conceptually}, to find the
algebraic equation for the Dyck paths considered by Volodia Retakh. To get the  analogous sequence
we only need to change the program by {\bf one line}. Let's call the analogous quantities $u_1(m,k)$ and $d_1(m,k)$.

{\tt u1:=proc(m,k) local r: option remember: if (k>1 and k mod 2=1) then  RETURN(0): fi: if m=0 then  0: else add(d1(m-r,k-r),r=1..k): fi: end:}

{\tt d1:=proc(m,k) local r: option remember:
if m=0 then  if k=0 then    RETURN(1):  else   RETURN(0): fi: fi: add(u1(m-r,k+r),r=1..m): end:}

In other words, just declaring that $u_1(m,k)=0$ if the elevation $k$ is an odd integer larger than $1$.

Typing

{\tt seq(d1(2*n,0),n=0..N);}

will let us get, very fast, the first $N+1$ terms, that would enable us to guess the algebraic equation satisfied by 
the generating function, that we can justify {\it a posteriori}, since we know that it {\bf exists}, saving us
the mental agony of doing it logically, by figuring out the intricate `grammar'.

{\bf The general case}

Since it is so  easy to tweak the  numerical dynamic programming procedure, why not be as general as can be?
Let $A$, $B$, $C$, $D$ be {\it arbitrary} sets of positive integers, either finite sets, or infinite sets
(like in Retakh's case) that are arithmetical progressions (or unions thereof). We are interested
in counting Dyck paths that obey the following restrictions

$\bullet$ No peak  can be of a  height that belongs to $A$.

$\bullet$ No valley  can be of a  height that belongs to $B$.

$\bullet$ No upward run  can be of  a length that belongs to $C$.

$\bullet$ No downward run  can be of a length that belongs to $D$.

Then we declare that $u(m,k)=0$ if $k \in A$ and $d(m,k)=0$ if $k \in B$ and otherwise
$$
d(m,k) \, = \,\sum_
{{1 \leq r \leq m}   \atop {r \not \in C}}
u(m-r,k+r) \quad .
$$
Analogously,
$$
u(m,k) \, = \, \sum_
{{1 \leq r \leq m}   \atop {r \not \in D}}
d(m-r,k-r) \quad .
$$

Then we get, very fast, sufficiently many terms to guess an algebraic equation, by  finding $\{d(2n,0)\}_{n=0}^{N}$,
for a sufficiently large $N$.

{\bf Guessing linear recurrences}

It is well-known (see [KP], Theorem 6.1) that if $f(x)$ is an algebraic formal power series
(like in our case), then it satisfies a linear differential equation with polynomial coefficients,
i.e. it is $D$-finite, and hence its sequence of coefficients, $\{a(n)\}_{n=0}^{\infty}$, satisfies a linear
recurrence equation with polynomial coefficients, i.e. is $P$-recursive.
While there are easy algorithms for finding these, they do not always give the minimal recurrence,
and once again, let's keep it simple! Just guess such a recurrence using {\it undetermined coefficients},
and we are guaranteed by the background  `general nonsense' that everything is rigorously proved,
and we don't have to worry about Richard Guy's Strong Law of Small Numbers.

\vfill\eject

{\bf The Maple package {\tt Dyck.txt}}

Everything here is implemented in the Maple package  {\tt Dyck.txt} available from the front of this
article

{\tt https://sites.math.rutgers.edu/\~{}zeilberg/mamarim/mamarimhtml/dyck.html} \quad .

There you would also find long web-books with many deep enumeration theorems. Let us present
just one {\bf random} example.

Typing 

{\tt Theorem($\{1\}$,$\{\}$,$\{2\}$,$\{1\}$,60,P,x,n,a,20,1000);} 

gives

{\bf Sample Theorem}: Let $a(n)$ be the number of Dyck paths of semi-length $n$ obeying the following restrictions.
The height of a peak can not belong to  $\{1\}$, no upward-run can have a length belonging to $\{2\}$, and
no downward-run can have a length that belongs to $\{1\}$, then the generating function
$$
f(x) \, := \,  \sum_{n=0}^{\infty} a(n) x^n \quad ,
$$
satisfies the algebraic equation
$$
1+  \left( -{x}^{2}-x-1 \right) f \left( x \right) +
\left( {x}^{4}+{x}^{3}+{x}^{2}+x \right)   f \left( x \right)^{2}  =0 \quad,
$$
and the sequence $a(n)$ satisfies the following linear recurrence
$$
a \left( n \right) ={\frac { \left( n-2 \right) a \left( n-1 \right) }{n+1}}+2\,{\frac { \left( n-2 \right) a \left( n-2 \right) }{n+1}}+{\frac { \left( 4\,
n-11 \right) a \left( n-3 \right) }{n+1}}
$$
$$
+{\frac { \left( 8\,n-25 \right) a \left( n-4 \right) }{n+1}}+6\,{\frac { \left( n-4 \right) a \left( n-5 \right) 
}{n+1}}+{\frac { \left( 5\,n-22 \right) a \left( n-6 \right) }{n+1}}+3\,{\frac { \left( n-5 \right) a \left( n-7 \right) }{n+1}} \quad,
$$
subject to the initial conditions $a(1) = 0, a(2) = 0, a(3) = 1, a(4) = 2, a(5) = 3, a(6) = 7, a(7) = 1 $.

{\bf Part II: Symbolic Dynamic Programming}

We will now briefly describe the {\bf clever} way, implemented in the Maple package {\tt DyckClever.txt}.
Note that the ``thinking" and ``research" is done by the computer {\it all by itself}, using
{\it symbolic} dynamic programming, to generate a {\bf system of algebraic equations}, and then
it solves that system. What happens is that in order to study the generating function for the
original set, we are {\bf forced} to consider other sets, that in turn, necessitate yet more sets.
Sooner or later, if all goes well, there would be no more new `uninvited guests', and 
the computer would have a {\bf finite} system of algebraic equations 
with {\bf as many equations as unknowns}. Using the Buchberger algorithm {\it under the hood}, it solves that system,
giving us much more than we asked for, not just the desired generating function, but
lots of other ones that we had to introduce, and that we may  not care about.

{\bf Avoiding Peak-Heights and Valley Heights with Finite Sets to Avoid}

Suppose that we have two {\bf arbitrary} {\bf finite} sets of non-negative  integers $A$ and $B$, and we are interested
in $f_{A,B}(x)$, the ordinary generating function enumerating Dyck paths such that

$\bullet$ None of the peak-heights is in $A$ .

$\bullet$ None of the  valley-heights is in $B$ .

Assume, for now, that $0 \not \in A$ and $0 \not \in B$, and define
$$
A_1 = \{ a-1 \, : \, a \in A \} \quad ,
$$
$$
B_1 = \{ b-1 \, : \, b \in B \} \quad .
$$

Recall that every  Dyck $P$  is either empty or else can be written as
$$
U P_1 D P_2 \quad,
$$
where $P_1$ and $P_2$ are Dyck  paths on their own right. If $P$ is counted by $f_{A,B}(x)$ then
$P_1$ is counted by $f_{A_1,B_1}(x)$, but $P_2$ is counted again by $f_{A,B}(x)$. This
leads to the quadratic equation
$$
f_{A,B}(x) \, = \,
1+ x f_{A_1,B_1}(x) f_{A,B}(x) \quad .
$$

Alas now we have to set-up an equation for $f_{A_1,B_1}(x)$ and keep going.
Sooner or later we will get an $f_{A',B'}(x)$ where either $A'$ or $B'$ contain $0$ (or both).
It is readily seen that if $0 \in A$ , then writing $A_1:=A \backslash \{0\}$, we get the
equation
$$
f_{A,B}(x) \, = \,f_{A_1,B}(x)-1 \quad .
$$

If $0 \not \in A$ but $0 \in B$ then define

$$
A_1 = \{ a-1 \, : \, a \in A \} \quad ,
$$
as above and
$$
B_1 = \{ b-1 \, : \, b \in B \} \backslash \{-1\} \quad ,
$$
and the corresponding equation is
$$
f_{A,B}(x)= 1+ x f_{A_1,B_1}(x)  \quad .
$$

Sooner or later there would be no more new `states' $[A,B]$ and we would have
a finite set of quadratic equations with as many equations as unknowns. Solving the resulting system
we would get, {\it inter-alia}, our original object of desire, $f_{A,B}(x)$.

This is implemented in  procedure {\tt fAB(A,B,x,P)} in the Maple package {\tt DyckClever.txt} available from

{\tt https://sites.math.rutgers.edu/\~{}zeilberg/tokhniot/DyckClever.txt} \quad .

Just to take a random example, in order to get the quadratic equation satisfied by
$P=\sum_{n=0}^{\infty} a(n)x^n$,  where $a(n)$ is the number of Dyck paths with
no peak-heights in the set $\{2,5\}$ and no valley-heights in the set $\{1,4\}$, type

{\tt fAB($\{$2,5$\}$,$\{1$,4$\}$,x,P);}

getting, in a fraction of a second,
$$
 \left( 4\,{x}^{6}-13\,{x}^{5}+24\,{x}^{4}-27\,{x}^{3}+19\,{x}^{2}-7\,x+1 \right) {P}^{2}
$$
$$
+ \left( 4\,{x}^{5}-15\,{x}^{4}+26\,{x}^{3}-26\,{x}^{2}+12\,x-2 \right) P+{x}^{4}-4\,{x}^{3}+8\,{x}^{2}-5\,x+1 \, = \, 0 \quad .
$$

The closely related procedure {\tt fABcat(A,B,x,C)} expresses $P$ in terms of the Catalan generating function
$C:=(1-\sqrt{1-4x})/(2x)$. Typing

{\tt fABcat( $\{$ 2,5  $\}$ , $\{$ 1,4 $\}$,x,C);}

would give
$$
{\frac {-C{x}^{2}+{x}^{2}-2\,x+1}{ \left( {x}^{3}-{x}^{2} \right) C-2\,{x}^{3}+3\,{x}^{2}-3\,x+1}} \quad .
$$

{\bf Avoiding Peak-Heights and Valley Heights with infinite Sets to Avoid}

The original motivating problem, asked by Volodia Retakh (see above) asked for the
number of Dyck paths whose peak-heights is either $1$ or even, in other words the
set of Dyck paths none of whose peak-heights is in the range of the arithmetical progression
$2r+3$ for $r \geq 0$. The above procedure $fAB(A,B,x,P)$ can be easily modified to
procedure

{\tt fABr(A,B,r,x,P)} \quad ,

where $A$ and $B$ are sets of arithmetical progressions written in the form $ar+b$ for $a$ and $b$ non-negative integers
and $r$ is a {\bf symbol} ranging over the non-negative integers. For example, to get 
Retakh's result type

{\tt fABr($\{$2*r+3 $\}$,$\{\}$,r,x,P);} \quad ,

getting
$$
{x}^{2}{P}^{2}+Px-P+1 \, = \, 0 \quad .
$$

For the equation satisfied by the generating function of the sequence enumerating Dyck paths
where neither peak-heights nor valley-height is  a member of the arithmetical progression $5r+1$  type

{\tt fABr( $\{$ 5*r+1 $\}$ , $\{$ 5*r+1 $\}$ ,r,x,P);}

getting
$$
\left( {x}^{3}+{x}^{2}-x \right) {P}^{2}+ \left( {x}^{3}-2\,{x}^{2}-x+1 \right) P+2\,x-1 \, = \, 0 \quad  .
$$

{\bf Avoiding Ascending Run-Lengths  and Descending Run-Lengths with Finite Sets to Avoid}

An {\bf irreducible} Dyck path of {\it semi-length} $n$ is one who never touches the $x$-axis except, of course, at
the beginning $(0,0$) and the end, when it is at $(2n,0)$.

Let $C$ and $D$ be {\bf arbitrary} finite sets of positive integers. We are
interested in finding the algebraic equation satisfied by the
generating function of the enumerating sequence of Dyck paths that
never have an {\bf ascending run-length} in $C$ and
never have a {\bf descending run-length} in $D$. Let's call it $h_{C,D}(x)$.

Alas, in order to get the {\it symbolic dynamic programming} (recursive) procedure  rolling, we
are {\bf forced} to consider more general classes.

Given sets of positive integers $C,C_1,D,D_1$,
let $h_{C,C_1,D,D_1}(x)$ be the weight-enumerators with $weight(P):=x^{SemiLength(P)}$ of {\bf all} Dyck paths $P$ such that

$\bullet$ The length of a starting upward-run is {\bf not} in $C_1$ .

$\bullet$ All other  upward run-lengths are {\bf not} in $C$.

$\bullet$ The length of an ending downward-run is {\bf not} in $D_1$.

$\bullet$ All other  downward run-lengths are {\bf not} in $D$.

Let $H_{C,C_1,D,D_1}(x)$ be the analogous quantity enumerating  {\bf irreducible} Dyck paths with these restrictions.

We have the following equation
$$
h_{C,C_1,D,D_1}(x) \, = \,
1 \, +\,  H_{C,C_1,D,D_1}(x)  \, + \,
 H_{C,C_1,D,D}(x)\,h_{C,C,D,D}(x)\,H_{C,C,D,D_1}(x) \quad .
$$
This follows from the fact that every such Dyck path is either the empty path (weight $1$),
or is irreducible, or else can be viewed as concatenation of three paths, 
such that the first is weight-counted  by  $H_{C,C_1,D,D}(x)$, the second (possibly empty) path counted
by $h_{C,C,D,D}(x)$ and the third path by $H_{C,C,D,D_1}(x)$.

We must now find an equation for $H_{C,C_1,D,D_1}$.

If $1 \not \in C_1$ and $1 \not \in D_1$ then, let
$$
C_2:= \{ c-1 \, : \, c \in C_1\} \quad ,
$$
$$
D_2:= \{ d-1 \, : \, d \in D_1\} \quad ,
$$

Obviously, we have
$$
H_{C,C_1,D,D_1}(x) \, = \, x \, h_{C,C_2,D,D_2}(x)  \quad .
$$

If $1 \in C_1$ and $1 \not \in D_1$ then
$$
H_{C,C_1,D,D_1} (x) \, = \,  H_{C,C_1 \backslash \{1\},D,D_1} (x) -x \quad ,
$$
since the set of Dyck paths counted by  $H_{C,C_1 \backslash \{1\},D,D_1} (x)$ and  $H_{C,C_1 ,D,D_1} (x)$
are almost the same, the only exception is $UD$ that belongs to the former, but not to the latter, and whose weight is $1$.

Similarly, if $1 \not \in C_1$ and $1 \in D_1$ then
$$
H_{C,C_1,D,D_1} (x) \, = \,  H_{C,C_1,D,D_1 \backslash \{1\}} (x) -x \quad ,
$$

Finally, if $1 \in C_1$ and $1 \in D_1$ then
$$
H_{C,C_1,D,D_1} (x) \, = \,  H_{C,C_1  \backslash \{1\}, D,D_1 \backslash \{1\} } (x) -x \quad .
$$

This is implemented in procedure {\tt fCD(C,D,x,P)}. 
Just to take a random example, in order to get the algebraic equation satisfied by
$P=\sum_{n=0}^{\infty} a(n)x^n$  where $a(n)$ is the number of Dyck paths with
no ascending run of length $2$  and no  descending run of length $3$ just type

{\tt fCD($\{$2$\}$,$\{$     3 $\}$,x,P);} 

getting
$$
1- \left( x+1 \right)  \left( {x}^{2}-x+1 \right) P+x \left( {x}^{3}+{x}^{2}-x+1 \right) {P}^{2}+{x}^{4} \left( {x}^{5}-{x}^{3}+2\,{x}^{2}-1 \right) {P}^{3}
$$
$$
+{x}^{4} \left( {x}^{6}-3\,{x}^{5}+{x}^{4}-x+1 \right) {P}^{4}+{x}^{6} \left( {x}^{7}+2\,{x}^{4}+{x}^{3}+{x}^{2}+1 \right) {P}^{5}
$$
$$
+{x}^{9} \left( 2\,{x}^{6}
+{x}^{5}+2\,{x}^{3}+3\,{x}^{2}-x+1 \right) {P}^{6}+{x}^{11} \left( {x}^{6}+2\,{x}^{5}+{x}^{3}+3\,{x}^{2}-2\,x+1 \right) {P}^{7}
$$
$$
+{x}^{15} \left( {x}^{4}+{x}^{3}+{x}^{2}+1 \right) {P}^{8}+{x}^{18} \left( {x}^{3}+{x}^{2}+1 \right) {P}^{9}+{x}^{22}{P}^{10} \, = \, 0 \quad .
$$

{\bf Avoiding Ascending Run-Lengths and Descending Run-Lengths with infinite Sets to Avoid}

The above procedure $fCD(C,D,x,P)$ can be easily modified to
procedure

{\tt fCDr(C,D,r,x,P)}

where $C$ and $D$ are sets of arithmetical progressions written in the form $ar+b$ for $a$ and $b$ non-negative integers
and $r$ is a symbol ranging over the non-negative integers. For example, to get  the generating function for
the sequence of Dyck paths where there are no ascending run-lengths in the infinite set $\{3,5,7,9,11, \dots\}$, type:

{\tt fCDr($\{$2*r+3 $\}$,$\{\}$,r,x,P);}

getting
$$
-1+ \left( -x+1 \right) P+{x}^{2} \left( x-1 \right) {P}^{3} \, = \, 0 \quad .
$$

For the equation satisfied by the generating function of the sequence enumerating Dyck paths
where there are no ascending run-lengths in the infinite set $\{2,4,6,8,10, \dots\}$, 
and no descending run-lengths in the infinite set $\{1,3,5,7,9, \dots\}$,  type:

{\tt fABr( $\{$ 2*r+2 $\}$ , $\{$ 2*r+1 $\}$ ,r,x,P);}

getting
$$
1+ \left( -{x}^{2}-1 \right) P-{P}^{2}{x}^{2}+{x}^{2} \left( 3\,{x}^{2}+2 \right) {P}^{3}-{P}^{4}{x}^{4}-{x}^{4} \left( {x}^{2}+1 \right) {P}^{5}
+ x^6 {P}^{6} \, = \, 0 \quad .
$$

{\bf Homework Assignments}

As in most of our papers, our main goal was to demonstrate an {\it approach} rather than discover new theorems.
It would be relatively straightforward to combine procedures {\tt fAB(A,B,x,P)} and
{\tt fCD(C,D,x,P)} to produce procedure  {\tt fABCD(A,B,C,D,x,P)} using the `clever' approach,
thereby giving a `clever analog' of the general problem discussed in Part I.
We leave this to the interested reader.

Avoiding an ascending run-length of length $a$ is the same thing as saying that $U^aD$ can't occur 
at the start of the word,  and $DU^aD$ can't  occur later, and analogously for avoiding descending  run-lengths.
This naturally generalizes to the problem of counting Dyck paths
avoiding as {\it consecutive subwords} (i.e. `factors' in the language of formal languages) 
an arbitrary  finite set of {\bf forbidden subwords}.
Both the `dumb' approach and the `clever' approach can be easily adapted to this problem,
and even when the subwords to avoid are `infinite set' consisting of regular expressions. 
We also leave this to the interested reader, as well as {\it interfacing} it with restricted peak-heights
and valley-heights.

{\bf Conclusion}

For many enumeration problems the `{\it dumb}' approach, {\it keeping it simple}, suffices, 
notwithstanding Richard Guy's cautionary tales, since by {\it mumbling} a few words, one can
easily  shut-up the traditional, machinophobic  `rigorist', who would have to accept it,
at least reluctantly (he would still say that such a proof `gives no insight').
This KISS approach is based on {\it number-crunching}, using {\it numeric} dynamic programming as the engine.

But for larger problems, we may need, unfortunately,  to do some `thinking' 
by examining the combinatorial-recursive structure of the combintorial sets involved.
Luckily, this approach can be automated as well, requiring us to only meta-think {\bf once} by
writing a {\bf general} program that {\it teaches} the computer 
to do the thinking for us.
This {\it clever} approach is based on {\it symbol-crunching}, using {\it symbolic} dynamic programming as the engine.

{\bf References}

[ELY] Sen-Peng Eu,  Shu-Chung Liu, and  Yeong-Nan Yeh, {\it Dyck Paths with Peaks Avoiding or Restricted to a Given Set},
Studies in Applied Mathematics {\bf 111} (2003), 453-465.

[G1] Richard K. Guy, {\it The strong  law of small numbers}, Amer. Math. Monthly, {\bf 95} (1988), 697--712.

[G2] Richard K. Guy, {\it The second strong  law of small numbers}, Math. Mag. {\bf 63} (1990), 3--20.

[KP] Manuel Kauers  and Peter Paule, {\it ``The Concrete Tetrahedron''}, Springer, 2011.

[M] Toufik Mansour, {\it Counting Peaks at Height k in a Dyck Path},     Journal of Integer Sequences, {\bf 5} (2002), Article 02.1.1

[PW] P. Peart and W.-J. Woan, {\it Dyck Paths With No Peaks at Height k}, J. Integer Sequences {\bf 4} (2001), \#01.1.3.

[Sl] Neil Sloane, {\it The On-Line Encyclopedia  of Integer Sequences},  {\tt https://oeis.org} \quad .

[St] Richard Stanley, {\it ``Catalan Numbers''}, Cambridge University Press, 2015.

[Z1] Doron Zeilberger, {\it Lagrange Inversion Without Tears (Analysis) (based on Henrici)},
The Personal Journal of Shalosh B. Ekhad and Doron Zeilberger, \hfill\break
{\tt https://sites.math.rutgers.edu/\~{}zeilberg/mamarim/mamarimhtml/lag.html} \quad .

\bigskip
\bigskip
\hrule
\bigskip
Shalosh B. Ekhad, c/o D. Zeilberger, Department of Mathematics, Rutgers University (New Brunswick), Hill Center-Busch Campus, 110 Frelinghuysen
Rd., Piscataway, NJ 08854-8019, USA. \hfill\break
Email: {\tt ShaloshBEkhad at gmail dot com}   \quad .
\bigskip
Doron Zeilberger, Department of Mathematics, Rutgers University (New Brunswick), Hill Center-Busch Campus, 110 Frelinghuysen
Rd., Piscataway, NJ 08854-8019, USA. \hfill\break
Email: {\tt DoronZeil at gmail  dot com}   \quad .
\bigskip
\hrule
\bigskip
{\bf Exclusively published in the Personal Journal of Shalosh B.  Ekhad and Doron Zeilberger and arxiv.org \quad .}
\bigskip
\hrule
\bigskip
First Written: June 3, 2020. This version: June 12, 2020.
\end